\documentclass[12pt,a4paper,fleqn]{article}

\usepackage{amssymb,amsmath,amsfonts}

\makeatletter

\@addtoreset{equation}{section}

\newtheorem{lem}{\sc Lemma }[section]
\newtheorem{theo}[lem]{\sc Theorem}
\newtheorem{propo}[lem]{\sc Proposition}
\newtheorem{coro}[lem]{\sc Corollary}

\newtheorem{exa}[lem]{\sc Example}

\newcommand{\qed}{\hspace*{2mm} \hfill $\Box $}

\newcommand{\vp}{\varphi} 

\newcommand{\ran}{\rangle}
\newcommand{\lan}{\langle}
\newcommand{\ra}{\rightarrow}

\newcommand{\R}{\mathbb R}

\newcommand{\ds}{\displaystyle}

\newcommand{\be}{\begin{equation}}
\newcommand{\ee}{\end{equation}}

\newcommand{\ov}{\overline}

\newcommand{\la}{\lambda} \newcommand{\La}{\Lambda}

\newcommand{\ba}{\begin{array}}
\newcommand{\ea}{\end{array}}

\newcommand{\bmt}{\left(\begin{array}}
\newcommand{\emt}{\end{array}\right)}

\title{A Note on Duality in  \\ Reverse Convex Optimization}

\author{Joachim Gwinner\thanks{
Institut  f\"ur Angewandte Mathematik, 
Fakult\"at  f\"ur  Luft- und Raumfahrttechnik,
Universit\"at der Bundeswehr M\"unchen, D--85577 Neubiberg/M\"unchen 
 {\it joachim.gwinner@unibw-muenchen.de}}}

\date{}

\begin{document}

\maketitle
 
\centerline{Dedicated to the memory of Prof. W. Oettli and Prof. D. Pallaschke}

\setcounter{page}{1}

\begin{abstract}
{ \small {\sc Abstract.}
In this note we explore duality in reverse convex optimization
with reverse convex inequality constraints.
While we are examining the special case of a finite index set of the
inequality constraints, we are primarily interested in the general case 
of an arbitrary infinite index set with neither extra topological structure nor 
extra measure structure.
 }
\end{abstract}


\begin{quote}
\end{quote}

 \section{Introduction}  

\setcounter{equation}{0}

This note is concerned with duality in reverse convex optimization problems of the form
\begin{eqnarray*}
(P_{(>,\ge)}) \quad
\left\{ \begin{array}{ll}
\textnormal{minimize } & f(x) \\[0.5ex]
\textnormal{subject to } & h(t,x) > 0 \textnormal{ for all } t \in T_{>}  \\[0.5ex]
\textnormal{and } & h(t,x) \ge  0 \textnormal{ for all } t \in T_{\ge}
\,,
 \end{array} \right.
\end{eqnarray*}
% DC f - g , subject to h - k LATER ?
where 
 $T= T_{>} \cup T_{\ge}   \not= \emptyset$ 
 with $ T_{>} \cap T_{\ge} = \emptyset$ is an index set.
and $f : X \ra \R,  h(t, \cdot): X \ra \R$  for all $t \in T$
are convex on a real topological vector space $X$. 
% on a convex set $X$ in a real linear space.
To include implicit constraints we can more generally admit $f$ to be extended real-valued, $f : X \ra \R \cup \{\infty\}$, which is proper, that is,
$\textnormal{dom }f = \{x \in X| f(x) < \infty \} \not= \emptyset$. 
While we are examining the special case of a finite index set $T$
with cardinality $|T| < \infty$, we are primarily interested in the general case  of an arbitrary infinite index set $T$ with no extra topological  or measure structure. Thus this note is close in spirit to the monograph \cite{PalRol} of Pallaschke and Rolewicz that develops the mathematical
foundations of optimization theory in a general and modern setting
dispensing with unnecessary structure.

Reverse convex programming can be traced back to \cite{HilJac}
and the references given therein.  See also the monographs of 
R. Horst and Hoang Tuy \cite{HorTuy} and of Hoang Tuy \cite{Tu-16}
that set reverse convex optimization in perspective to global
optimization. A  duality theory for optimization problems with reverse
convex constraint sets has been developed in
 \cite{TBO,Sin-92,Str,Sin-00,Pen} in chronological order.
%  see also \cite{GCM} that sets \cite{TBO} in perspective to more recent work % in nonconvex optimization.
\cite{Jey} provides dual characterizations of the containment of a closed convex set, defined by infinite convex constraints,  in a reverse-convex set,
defined by convex constraints, where these  dual characterizations are given in terms of epigraphs of conjugate functions.
\cite{AbAd} presents optimality conditions for a class of nonsmooth generalized semiinfinite programming problems  whose feasible sets are
represented by a single reverse convex constraint; that is, by a constraint of the form $h(x) \ge 0$ with a real-valued convex function $h$.
\cite{Tha} studies convex programs with several additional reverse convex constraints. 
While \cite{Le-1998,LemVol} establish duality results for 
nonconvex minimization problems with a single inequality constraint,  
\cite{MLVol} extend \cite{Le-1998,LemVol} to D.C. programming with several D.C. constraints.

First in this note we explore under which conditions the optimization
problems $(P_{(>)})$ and $(P_{(\ge)})$ are equivalent,
that is, their optimal values are equal, see Proposition \ref{prop2.1},
which is an preliminary result for our subsequent duality theory. Moreover,
Example \ref{ex2.2}  shows how a class of set reverse convex optimization problems can be recast as a reverse convex optimization problem of the form $(P_\ge)$ with a single inequality constraint and moreover,
in virtue of the Proposition \ref{prop2.1} also as a reverse convex optimization problem of the form $(P_>)$ with a single strict inequality constraint. Then in the main body of this note we treat the optimization problem of minimizing a convex function subject to an infinite
number of reverse convex constraints. To achieve a duality result for such a general class of nonconvex problems, see Theorem \ref{theo:3.2}, we employ regularization similar to \cite{GwiOet-89}. 
Furthermore in the case of a finite number of  reverse convex constraints,
Corollary \ref{coro:3.3} recovers a duality result due to Mart\'{\i}nez-Legaz, J.-E. and Volle \cite{MLVol}.

The key to the treatment of an infinite number of constraints is the embedding 
in the function space $\R^T$ endowed with the product topology which goes back to Aubin \cite{Aub} in his approach to minimax theorems. This approach  has been then further exploited  in  \cite{JeyGwi} on
solvability theorems for general inequality
systems applied to constrained optimization problems with an infinite number of constraints and to minimax problems on noncompact sets
and in \cite{GwiOet-94} that presents general  theorems of the alternative and extends duality theory of convex mathematical programming to inf-sup problems
under an abstract closedness assumption, thus avoiding the usual compactness requirement.

The outline of the note is as follows. The subsequent section 2 provides
Proposition \ref{prop2.1} relating $(P_{(>)})$ and $(P_{(\ge)})$,
moreover Example \ref{ex2.2}. The final section 3 presents the duality results 
Theorem \ref{theo:3.2} and  Corollary \ref{coro:3.3}.

\section{A relation between $(P_{(>)})$ and $(P_{(\ge)})$}

In this section we explore under which conditions the optimization
problems $(P_{(>)})$ and $(P_{(\ge)})$ are equivalent,
that is, there holds 
$\textnormal{inf} (P_{(>)}) = \textnormal{inf} (P_{(\ge)}) $
 for the optimal values 
$$
\textnormal{inf} (P_{(>)}) = \textnormal{inf }
\{f(x) | x \in X, h(t,x) > 0 \, \forall  t \in T \}
$$
and  
$$
\textnormal{inf} (P_{(\ge)}) = \textnormal{inf}
\{f(x) | x \in X, h(t,x) \ge 0 \, \forall  t \in T \}
$$
Here we do not need a topological structure on $X$, so $X$ is more generally a linear space or a convex set thereof. 

\begin{propo} \label{prop2.1}
Assume inf$(P_{(\ge)}) > - \infty$.

(i) Then the optimal values are equal, if and only if
for any $\eta > 0$  there exists $ x_\eta \in X$ such that 
$\textnormal{inf} (P_{(\ge)}) > f(x_\eta) - \eta$ 
and $h(t,x_\eta) > 0 \,\, \forall  t \in T$.

(ii) Suppose
\begin{eqnarray*} 
&& T \textnormal{ is finite}, X  \textnormal{ is a linear space}, \\ 
&& f \textnormal{ is upper semicontinuous on line segments}, \\ 
&& h(t,\cdot)  \textnormal{ is lower semicontinuous on line segments} 
 \, \forall  t \in T \,, \\ 
&& \textnormal{and there exists } \tilde x \in X  \textnormal{ such that } 
 h(t,\tilde x) < 0 \, \forall  t \in T \,.
\end{eqnarray*}
 Then the optimal values 
$\textnormal{inf} (P_{(>)})$ and 
$\textnormal{inf} (P_{(\ge)})$ coincide.
\end{propo}

{\sc Proof}. Obviously inf$(P_{(\ge)}) \le$ inf$(P_{(>)})$. Let 
$\delta := \textnormal{inf} (P_{(>)}) - \textnormal{inf} (P_{(\ge)}).$

(i) To prove the sufficiency, assume $\delta >0$. Choose $\eta \in (0, \delta/2)$. Then $x_\eta$ is feasible for
 $(P_{(>)})$, and so inf$(P_{(>)}) \le f(x_\eta)$. Hence 
\begin{eqnarray*}
 &\textnormal{inf} (P_{(\ge)})&  > f(x_\eta) - \eta \ge
 \textnormal{inf} (P_{(>)}) - \eta \\
&& > \textnormal{inf} (P_{(>)}) - \ds \frac{1}{2} \, \textnormal{inf}(P_{(>)})
+  \ds \frac{1}{2} \, \textnormal{inf} (P_{(\ge}) \\
&& =  \ds \frac{1}{2} \, [\textnormal{inf} (P_{(>)}) + \textnormal{inf} (P_{(\ge)}]\,,
\end{eqnarray*}
what leads to $\textnormal{inf} (P_{(\ge)}) > \textnormal{inf} (P_{(>)}) $,
what is absurd.

The necessity is trivial.

(ii) Again assume $\delta >0$ Since inf$(P_{(\ge)}) > - \infty$, there exists $\ov x \in X$ 
with $ h(t,\ov x) \ge 0 \,\, \forall  t \in T$ such that
$f(\ov x) < \textnormal{inf} (P_{(\ge)}) - \ds \frac{1}{2} \delta $.
If $ h(t,\ov x) > 0 \,\, \forall  t \in T$, then the above argument in $(i)$ 
leads to a contradiction again. Otherwise 
${\ov T}_= := \{ t  \in T| h(t, \ov x) = 0 \}$ is not empty;
 let $\ov T_> := T \setminus  \ov T_=$, 
$x_s := \ov x + s ( \ov x - \tilde x) \in X$ for $s > 0$. 
Then $\ov x =  \frac{1}{1+s} x_s +  \frac{s}{1+s} \tilde x $, hence by convexity of $h(t,\cdot)$, for all $t \in T$,
$$
h(t,\ov x) \le \ds \frac{1}{1+s} h(t,x_s) + \ds\frac{s}{1+s} h(t,\tilde x) \,.
$$
This implies for any $s>0$, for all $t \in \ov T_=$, $h(t,x_s) > 0$;
since otherwise $h(t,x_s) \le 0$ would result in
$h(t,\ov x) \le \frac{s}{1+s} h(t,\tilde x) < 0$, a contradiction.

Next, since $h(t,\cdot)$ is lower semicontinuous on line segments
and $\ov T_>$ is finite, there exists $s_1 > 0$ such that
$h(t,x_s) > 0$ for $0 < s < s_1$ and for all $t \in \ov T_>$. 

Further, since $f$ is upper semicontinuous on line segments,  
there exists $s_2 > 0$ such that
$f(x_s)  < \textnormal{inf} (P_{(\ge)}) - \ds \frac{1}{2} \delta $
 for $0 < s < s_2$. Let $\hat s := \ds \frac{1}{2} \min(s_1,s_2),
\hat x_s :=  \ov x + \hat s ( \ov x - \tilde x)$.
Then we obtain 
$f(\hat x_s)  < \textnormal{inf} (P_{(\ge)}) - \ds \frac{1}{2} \delta $
and $h(t,\hat x_s) > 0$ for all $t \in T$, thus a contradiction as in case (i). 
 \qed

In the following example we show how a class of set reverse convex optimization problems can be recast as a reverse convex optimization problem of the form $(P_\ge)$ with a single inequality constraint and moreover,
in virtue of the above Proposition \ref{prop2.1} also as a reverse convex optimization problem of the form $(P_>)$ with a single strict inequality constraint.

\begin{exa} \label{ex2.2}
Consider the optimization problem 
  
$$
  (Q)\quad
 \min f(x) \mbox{ such that } T x \not \in \mbox{int } D,
  $$

where  
  $$
  \begin{array}{l}
    f: \R^n \rightarrow \R \cup \{+\infty\} \mbox{ lower semicontinuous
       (proper and convex)} \\
    T: \R^n \rightarrow \R^d \mbox{ linear} \\
    D \not = \emptyset \mbox{ closed convex } \subset \R^d \mbox{ with int } D
    \not= \emptyset ~.
  \end{array}
  $$

Assume that the set reverse convex constraint
$T x \not \in \mbox{int }D$ is solvable and there exists
 $\tilde x \in$ argmin $\{ f (x) : x \in \R^n \}$. 

We can focus to the case that $T \tilde x \in$ int $D$, since otherwise $\tilde x$ solves $(Q)$ and $(Q)$ is only a standard convex optimization problem. 
Let $V: = D - T \tilde x$. Then $V$ is convex with $0 \in$ int $D$ and the polar 
$E:= V^0$ is $\not = \emptyset$, convex and compact.

Now rewrite the reverse convex constraint by the separation theorem:
$T x \not \in \mbox{int } D 
\Leftrightarrow \; \exists \; u \in E: \langle u, T x - T \tilde x) \rangle \geq 1$.
Thus 
 \begin{eqnarray*}
   \inf (Q) & = &
     \inf \big\{ f(x) \; |  \; \exists \; u \in E: \langle u, 
      T x - T \tilde x \rangle \geq 1 \big\} \\
   & = & \inf \big\{ f(x) \, | \; h(x) \ge 0  \big\} 
    \,,
\end{eqnarray*}
where 

$$
h(x) := \max \big\{  \langle u, Tx - T \tilde x \rangle |
\, u \in E \big\}
\,- \, 1  $$
is convex and continuous. Next $h(\tilde x) = -1 < 0$. 
Now assume in addition that $f$ is upper semicontinuous on line segments.
Then in virtue of Proposition \ref{prop2.1},

$$\inf (Q) =  \inf  \big\{ f(x) \; | \; h(x) > 0  \big\} \,.$$ 

\end{exa}

\section{Duality results for $(P_{(>)})$ without and with regularity assumption}

in this section we study the optimization problem  $(P_{(>)})$; first we provide a feasibility lemma and then show duality results. We work in the setting of a Hausdorff locally convex real linear space $X$ which is paired 
with its topological dual $Y$, where $\lan \cdot,\cdot \ran$ 
denotes the  canonical pairing between $X$ and $Y$. We emphasize that 
$T$ is an arbitrary nonvoid index set. Since $T$ does not bear a structure,
we work wih the function space $\R^T$ endowed with the product topology $\pi$.
Note that $\pi$ is the topology of pointwise convergence 
and there is a local neighborhood base of each point of $\R^T$
in which each set is determined by a finite number of elements of $T$.   
The nonnegative dual cone to $\R^T$ is then given by

$$
\La = \{\la: T \ra \R_+ | \exists \tilde T \textnormal{ finite} \subset T:
\la_t \equiv \la(t) = 0 \,(\forall t \in T \setminus \tilde T  \} \,.
$$

For notational simplicity set $h_t \equiv h(t,\cdot)$. The feasible set of 
 $(P_{(>)})$ will be denoted by 
$$
{\cal F}(P_{(>)}) = \bigcap_{t \in T} \{x \in X| h_t(x) > 0 \} 
$$ 
By Fenchel conjugation, see e.g. \cite{Lau,PalRol}, 
$h_t = h_t^{**}$, hence
$$
h_t(x) > 0 \Leftrightarrow
h_t^{**}(x) = \sup_{y \in Y}[\lan y,x \ran - h_t^{*}(y)] > 0 
\Leftrightarrow \exists y \in Y: \lan y,x \ran - h_t^{*}(y) > 0 \,.
$$
This means for $x \in {\cal F}(P_{(>)})$ that for any $t \in T$
there exists some $y \in Y$ such that the above strict inequality holds.
In virtue of the axiom of choice, see e.g. \cite{Dieu}, there exists a 
mapping $\vp: T \ra Y$ such that for all $t \in T$,

$$
\lan \vp_t,x \ran \equiv \lan \vp(t), x \ran  > 
h_t^*(\vp_t) \equiv h_t^*(\vp(t)) \equiv h^*(t,\vp(t)) \,.
$$

Thus we are led to introduce the set

$$
\Phi := \Phi(f,h) := \{\vp \in Y^T | \exists \ov x \in 
\textnormal{ dom } f: \lan \vp_t, \ov x \ran > h_t^*(\vp_t) \, 
(\forall t \in T) \}
$$

and obtain the feasilbility lemma:

\begin{lem}

${\cal F}(P_{(>)}) \not= \emptyset \Leftrightarrow 
 \Phi(f,h) \not= \emptyset \,.$

\end{lem}

Next we set

$$
j(\vp) := \inf \{f(x) | \lan \vp_t,x \ran > h_t^*(\vp_t) 
(\forall t \in T) \} \,.
$$

Then we have for the primal value of $P_{(>)}$

$$
\alpha := \inf(P_{(>)}) = \inf \{j(\vp) | \vp \in \Phi \}
$$

Moreover we set

$$
\ov j(\vp) := \inf \{f(x) | \lan \vp_t,x \ran \ge h_t^*(\vp_t) 
(\forall t \in T) \} \,.
$$

Then clearly $\alpha \ge \ov \alpha$, where

$$
\ov \alpha : =  \inf \{\ov j(\vp) | \vp \in \Phi \} \,,
$$
what we call the regularized primal value of $P_{(>)}$ similar to \cite{GwiOet-89}. 
To come up with the associated Lagrangian, we observe for any 
$\vp \in \Phi$,

$$
\lan \vp_t,x \ran \ge h_t^*(\vp_t) \,\forall t \in T
\Leftrightarrow
\sum_{t \in T} \la_t [h_t^*(\vp_t) - \lan \vp_t,x \ran] \le 0 \, 
\forall \la \in \La \,.
$$

Hence 

$$
\ov \alpha = \inf_{\vp \in \Phi} \inf_{x \in X} \sup_{\la \in \La}
[f(x) + H(\vp; x,\la)] \,,
$$
where 
$$
 H(\vp; x,\la) := \sum_{t \in T} \la_t [h_t^*(\vp_t) - \lan \vp_t,x \ran] \,.
$$

To apply the general duality theory of \cite{GwiOet-89} we introduce
for any $\vp \in \Phi$
\begin{eqnarray*}
& C(\vp) := \{v \in \R^X | \exists \la \in \La : H(\vp;\cdot,\la) \ge v 
\textnormal{ on } X \} \,, \\
& D(\vp) := \{u \in \R^\La | \exists x \in X : H(\vp;x,\cdot) \le u 
\textnormal{ on } \La \} \,.
\end{eqnarray*}

Then obviously $C(\vp)$ and $D(\vp)$ are convex, however $C(\vp)$ is not $\pi-$closed in $\R^X$, in general. Therefore as in \cite{GwiOet-89}, we are led to introduce the regularized dual value

$$
\ov \beta := \inf_{\vp \in \Phi} \, \, \sup_{v \in \pi -cl \, C(\vp)} 
\, \, \inf_{x \in X}
[f(x) + v(x)] 
$$ 

in addition to the dual value

\begin{eqnarray*}
& \beta & := \inf_{\vp \in \Phi} \, \sup_{v \in C(\vp)} \, \inf_{x \in X}
[f(x) + v(x)] \\
&& = \inf_{\vp \in \Phi} \, \sup_{\la \in \La} \, \inf_{x \in X}
[f(x) + H(\vp;x,\la)] \\
&& = \inf_{\vp \in \Phi} \, \sup_{\la \in \La} \, \inf_{x \in X}
[f(x) - \sum_{t \in T} \la_t \lan \vp_t,x \ran 
+  \sum_{t \in T} \la_t h_t^*(\vp_t)] \,. 
\end{eqnarray*}  

Since

$$
\inf_{x \in X} [f(x) - \sum_{t \in T} \la_t \lan \vp_t,x \ran] 
= -f^*(\sum_{t \in T} \la_t \vp_t ) \,,
$$ 

finally 

$$
\beta =  \inf_{\vp \in \Phi} \, \sup_{\la \in \La}
[ \sum_{t \in T} \la_t h_t^*(\vp_t) - f^*(\sum_{t \in T} \la_t \vp_t ) ] \,.
$$
Thus the duality theory of \cite{GwiOet-89},\cite[Theorem 3]{GwiOet-94}
applies to conclude the duality theorem without 
any regularity assumption (also called constraint qualification)

\begin{theo}\label{theo:3.2}
There holds
$$
\alpha \le \ov \alpha = \ov \beta \le \beta \,,
$$
where equality trivially holds when $\alpha = \infty$ or $\beta = - \infty$.
\end{theo}

Next we are dealing with the question, when (except the trivial cases above)
equality $ \alpha = \beta$ with dual attainment does hold.
To obtain this equality we are going to show $\ov \alpha = \alpha$
and $\ov \beta = \beta$. 

Note that $H(\vp;x,\la)$ is linear with respect to $\la$ on the convex cone $\La$. Hence we are in he setting of partially homogeneous programs
in \cite{GwiOet-94}. Therefore to show $\ov \beta = \beta$ with dual attainment we can apply \cite[Corollary 2]{GwiOet-94} that demands
two conditions, firstly the Karlin-Eisenberg regularity assumption

$$
(RA1) \quad \la \in \La, H(\vp; x,\la) \ge 0 \, (\forall x \in X)
\Rightarrow \la  = 0 
$$
for any $\vp \in \Phi$.
Now by construction, for any $\vp \in \Phi$, there exists 
$\ov x \in \textnormal{ dom } f$ such that 
$\lan \vp_t, \ov x \ran > h_t^*(\vp_t) \, 
(\forall t \in T) \}$. 
Hence  
$$
 H(\vp; \ov x,\la) = \sum_{t \in T} \la_t [h_t^*(\vp_t) - \lan \vp_t,\ov x \ran] \ge 0 \Rightarrow \la_t = 0 \, (\forall t \in T) 
$$
and $(RA1)$ is satisfied.

The second condition of \cite[Corollary 2]{GwiOet-94} demands a compact
"base" of $\La$, that is a   
subset $B \subset  \La$ such that $0 \notin B$ and $\La = \R_+ \, B$.
Consider $B := \{ \la \in \La| \sum\la_t =1\}$. Then $0 \notin B$, 
however, $B$ is not compact, in general. This requires $T$ to be finite.
 
Next to obtain $\ov \alpha = \alpha$, we apply Proposition \ref{prop2.1}
to
$$
\tilde h(t,x) := \lan \vp_t, x\ran - h_t^*(\vp_t)
$$
for fixed $\vp \in \Phi$.

To sum up, for finite index set $T$ with cardinality $|T| =m$,
$\La = \R^m_+$, and 
$$
\Phi = \{ (y_1, \dots, y_m) \in Y^m |\exists \ov x \in 
\textnormal{ dom } g: \lan y_j, \ov x \ran > h_j^*(y_j) \, (j=1,\ldots,m) \}
$$
we recover a duality result of Mart\'{\i}nez-Legaz and Volle
 \cite[Theorem 5.1]{MLVol} (with $h= g_j = 0 \, (j=1,\ldots,m))$
 
\begin{coro}\label{coro:3.3}
 There holds under the continuity and convexity assumptions above
\begin{eqnarray*}
& \alpha & = \inf \{f(x) |x \in X, h_j(x) > 0  \, (j=1,\ldots,m)\}  \\
&& = \inf_{y \in \Phi} \, \max_{\la \in \R^m_+} 
[\sum_{j=1}^m \la_j \, h_j^*(y_j) 
- f^*(\sum_{j=1}^m \la_j y_j)] \,. 
\end{eqnarray*}  
\end{coro}

\section*{Disclosure statement}

  The author reports that there are no competing interests to declare.

\end{document}